\documentclass[12pt,a4paper]{article}
\usepackage{amssymb,amsmath}
\newtheorem{theorem}{Theorem}[section]
\newtheorem{corollary}[theorem]{Corollary}
\newtheorem{lemma}[theorem]{Lemma}
\newtheorem{proposition}[theorem]{Proposition}
\newtheorem{definition}[theorem]{Definition}
\newtheorem{remark}[theorem]{Remark}
\usepackage[T1]{fontenc}
\title{Jacobi fields and the stability of minimal foliations of  arbitrary codimension}
\author{Krzysztof Andrzejewski}
\date{}
\begin{document}
\def\CC{\mathbb C}
\def\RR{\mathbb R}
\def\calf{\mathcal F}
\def\ric{\operatorname{Ric}}
\def\tr{\operatorname{Tr}}
\newcommand\brac[1]{\langle{#1}|}
\newcommand\ket[1]{|{#1}\rangle}
\newcommand\bracket[2]{\langle{#1},{#2}\rangle}
\newcommand\codim{\operatorname{codim}}
\newcommand\divergence{\operatorname{div}}
\newcommand\vol{\operatorname{vol}}

\maketitle
\begin{abstract}
In this article, we investigate the  stability of leaves of minimal foliations of arbitrary codimension.
We also study relations between Jacobi fields and vector fields which preserves a  foliation
and we use  these results  to Killing  fields. 
\par
{\bf 2000 Mathematics Subject Classification:} 53C12,53C42.
\par
{\bf Key words and phrases:} minimal  foliations, mean curvature, Jacobi field, Killing field.
\end{abstract}
\section*{Introduction}
\label{s:1}
A leaf of a minimal foliation $\calf$   is stable if the second derivative of the volume functional
with respect to any compactly supported normal variational field is nonnegative.
In this article, by direct calculations, we show that any leaf of a minimal foliation  with an integrable orthogonal distribution is stable.
Next, we introduce Jacobi fields of foliations, and  we investigate  relations between
these fields and vector fields preserving a foliation (Propositions \ref{p1} and \ref{p2}). Using these relations,
we show directly (not using the notion of calibration)  that any Killing vector field preserves
a minimal foliation $\calf$ having all  leaves compact and an
integrable orthogonal distribution  (Corollary \ref{c2}).
We also show that a Killing field preserves two orthogonal  complementary minimal foliations on
a closed manifold (Corollary \ref{c3}).  Finally, we give some consequences of these results.
Throughout the paper everything (manifolds, distribution, metrics, etc.) is assumed to be
$C^{\infty}$-differentiable and oriented.
\section{Stability results}
Let $M$ be an $m$-dimensional oriented, connected Riemannian  manifold.
On $M$, we consider  a foliation $\calf$, and let $n=\dim \calf$. Let  $D$ denote the distribution
corresponding to $\calf$, i.e., $D=T\calf$,  and $D^{\bot}$  the distribution  which
is the  orthogonal complement of $D$, $l=\dim D^{\bot}=m-n$. We assume that they are orientable and
transversally orientable. Let $\langle \cdot ,\cdot\rangle$  represent a metric on $M$
and $\nabla$ denote the Levi-Civita connection of the metric. Let $\Gamma(D)$ and  $\nabla^{\top}$
denote the set of  all  vector fields tangent to $D$ and  the induced connection in $D$, respectively. Similarly,
we have $\Gamma(D^{\bot})$ and $\nabla^{\bot}$. Moreover, $L(\Gamma(D), \Gamma(D^\bot))$
denotes the set of all $C^\infty$-linear transformations with the induced inner product (see equation (\ref{e0})). 
\par Throughout this  paper, we will use the following index convention  $1\leq i,j,\ldots\leq n,$
\quad $n+1\leq\alpha,\beta,\ldots \leq m$. Repeated indices denote summation over their range.
Let us take a local orthonormal frame $\{e_1,\ldots,e_m\}$ adapted to $D,D^{\bot}$, i.e.,  $\{e_i\}$
are tangent  to $D$  and $\{e_{\alpha}\}$ are tangent to $D^{\bot}$.
Moreover, $\{e_1,\ldots,e_m\}$,$\{e_i\}$ and $\{e_\alpha\}$ are compatible with the
orientation of $M,D$  and $D^{\bot}$, respectively. Then,  for  $A,B\in L(\Gamma(D), \Gamma(D^\bot))$, we obtain
\begin{equation}
\label{e0}
\langle A,B\rangle=\langle A(e_i),B(e_i)\rangle.
\end{equation}
Finally,  if $v$ is a vector tangent to $M$, then we write $v=v^\top+v^\bot$, where 
$v^\top$  belongs to $D$ and $v^\bot$  to $D^\bot$. 
\par Define the shape operator  $A^V\in L(\Gamma(D),\Gamma(D))$ of
$\calf$ with respect to $V\in\Gamma(D^\bot)$ by
\[
A^V(X)=-(\nabla_XV)^\top \quad {\rm for }\quad X\in\Gamma(D).
\]
Then, using the notation  $A^\alpha=A^{e_\alpha}$,
we have that  the mean curvature vector field $H$ of $\calf$ is given by
\[
H=\tr (A^\alpha )e_\alpha.
\]
We say that $\calf$ is minimal if $H=0$, i.e., if each leaf of $\calf$ is a minimal submanifold of $M$.
For $V\in\Gamma(D^\bot)$, we define mappings $\alpha_V,\nabla^{\bot}V\in L(\Gamma(D), \Gamma(D^\bot))$ by
\[ \nabla^{\bot}V(X)= \nabla^{\bot}_XV, \quad\alpha_V(X)=[V,X]^\bot
\quad {\rm for}\quad X\in \Gamma( D),\]
and the
 field $R(V)=(R(e_i,V)e_i)^{\bot}$. Here $R$ denotes the curvature tensor of $M$.
\par
In  the next part of this article, we will need the following lemma.
\begin{lemma}[\cite{aw}]
\label{l1}
Let  $p\in M$, and  let $ \{e_1,\ldots,e_m\}$  be a local orthonormal frame field adapted to $D,D^\bot$, such that
$(\nabla_Xe_i)^{\top}(p)=0$ and $(\nabla_Xe_\alpha)^{\bot}(p)=0$ for any vector field $X$ on
$M$. Then we have  at the point $p$
\begin{align*}
e_\alpha({A^{\beta}}^i_j)=&(A^\beta A^\alpha)_j^i-\langle R(e_j,e_\alpha)e_i,e_\beta\rangle
\\
+&\langle(\nabla _{e_\alpha} e_\gamma)^{\top},e_j\rangle\langle e_i,(\nabla_{e_\gamma}e_\beta)^\top
\rangle-\langle \nabla_{e_j}(\nabla _{e_\alpha}e_\beta)^{\top},e_i\rangle.
\end{align*}
\end{lemma}
Under the notation of Lemma \ref{l1},  we have
\begin{corollary}
\label{c1}
If $\calf$ is a minimal foliation and $D^\bot$ is  integrable, then we have at the
point $p$
\[
-\tr (A^\alpha A^\beta) +\langle R(e_\alpha),e_\beta\rangle=
\langle(\nabla _{e_\alpha} e_\gamma)^{\top},(\nabla_{e_\beta}e_\gamma)^\top
\rangle-\divergence_L ((\nabla _{e_\alpha}e_\beta)^{\top});
\]
where $\divergence_L(X)=\langle \nabla_{e_i} X,e_i\rangle$, for $X\in\Gamma(D)$.
\end{corollary}
Now, for $V,W\in\Gamma(D^\bot)$, we introduce   an auxiliary function  by
\[f_{V,W}=\langle\nabla^{\bot}V,\nabla^{\bot}W\rangle+\langle R(V),W\rangle
-\langle A^V,A^W\rangle.\]
Then,  we have the following lemma.
\begin{lemma}
\label{l2}
If $\calf$ is a minimal foliation and $D^\bot$ is  integrable,  then
\[
f_{V,W}=\langle \alpha_V,\alpha_W\rangle-\divergence_L((\nabla_VW)^\top).
\]
\end{lemma}
{\it Proof.}
Take a basis as in Lemma \ref{l1}. Then  at the point  $p$ we have the equalities
\begin{align*}
&\langle\nabla^{\bot}V,\nabla^{\bot}W\rangle=\langle(\nabla_{e_i}V)^\bot,(\nabla_{e_i}W)^\bot\rangle\\
&=\langle e_i(V^\alpha)e_\alpha+V^\alpha(\nabla_{e_i}e_\alpha)^\bot,
e_i(W^\beta)e_\beta+W^\beta(\nabla_{e_i}e_\beta)^\bot\rangle\\
&=e_i(V^\alpha)e_i(W^\alpha),
\end{align*}
and
\begin{align*}
&\langle A^V,A^W\rangle=\langle (\nabla_{e_i}V)^{\top},(\nabla_{e_i}W)^{\top}\rangle\\
&=V^\alpha W^\beta\langle A^\alpha(e_i),A^\beta(e_i)\rangle=V^\alpha W^\beta
\tr(A^\alpha A^\beta),
\end{align*}
and also
\[
\langle R(V),W\rangle=V^\alpha W^\beta \langle R(e_\alpha) e_\beta\rangle.
\]
Thus we have
\begin{equation}
\label{e1}
f_{V,W}(p)=e_i(V^\alpha)e_i(W^\alpha)+V^\alpha W^\beta\langle R(e_\alpha),e_\beta\rangle-
V^\alpha W^\beta \tr (A^\alpha A^\beta)|_p.
\end{equation}
On the other hand, from Corollary \ref{c1},  we have at the point $p$
\begin{align*}
&-V^\alpha W^\beta \tr (A^\alpha A^\beta) +V^\alpha W^\beta
\langle R(e_\alpha),e_\beta\rangle\\
&=\langle(\nabla _{V} e_\gamma)^{\top},(\nabla_{W}e_\gamma)^\top
\rangle-V^\alpha W^\beta\divergence_L ((\nabla _{e_\alpha}e_\beta)^{\top})\\
&=\langle(\nabla _{V} e_\alpha)^{\top},(\nabla_{W}e_\alpha)^\top
\rangle-\divergence_L ((\nabla _{V}W)^{\top})\\
&+W^\beta(\nabla_{e_\alpha}e_\beta)^\top( V^\alpha)+V^\alpha(\nabla_{e_\alpha}e_\beta)^\top( W^\beta)\\
&=\langle(\nabla _{V} e_\alpha)^{\top},(\nabla_{W}e_\alpha)^\top
\rangle-\divergence_L ((\nabla_{V}W)^\top)\\
&+(\nabla_{W}e_\alpha)^\top( V^\alpha)+(\nabla_{V}e_\alpha)^\top( W^\alpha).
\end{align*}
Using  this and  (\ref{e1}),  we obtain
\begin{align*}
&f_{V,W}(p)=e_i(V^\alpha)e_i(W^\alpha)+
\langle(\nabla _{V} e_\alpha)^{\top},(\nabla_{W}e_\alpha)^\top
\rangle\\
&-\divergence_L ((\nabla _{V}W)^{\top})
+(\nabla_{W}e_\alpha)^\top( V^\alpha)+(\nabla_{V}e_\alpha)^\top( W^\alpha)|_p.
\end{align*}
Moreover, we have at the point $p$
\begin{align*}
&\langle \alpha_V,\alpha_W\rangle=\langle(\nabla_Ve_i-\nabla_{e_i}V)^\bot,
(\nabla_We_i-\nabla_{e_i}W)^\bot\rangle\\
&=(\langle \nabla_V e_\alpha,e_i\rangle+\langle \nabla_{e_i}V,e_\alpha\rangle)
(\langle \nabla_W e_\alpha,e_i\rangle+\langle \nabla_{e_i}W,e_\alpha\rangle)\\
&=(\langle\nabla_V e_\alpha,e_i\rangle+e_i(V^\alpha))
(\langle\nabla_W e_\alpha,e_i\rangle+e_i(W^\alpha)).
\end{align*}
Thus
\[
f_{V,W}(p)=\langle\alpha_V,\alpha_W\rangle(p)-\divergence_L((\nabla_VW)^\top)(p).
\]
Since $p$ is arbitrary, we complete  the proof.\hfill$\square$
\par Now, let $L$ be a leaf of the foliation $\calf$ with the induced metric and the Levi-Civita connection
$\tilde\nabla$. Similarly as before, we introduce the connection in $\Gamma((TL)^\bot)$, the second
fundamental form $\tilde A^{\tilde V}$ of a leaf $L$ and $\tilde\nabla^\bot\tilde V$
for an arbitrary $\tilde V\in\Gamma((TL)^\bot)$ (see \cite{s1}). Let $L$ be a minimal submanifold
of $M$, we say that $L$ is {\it stable}  if the  inequality
\[
\int_L f_{\tilde V}\geq 0
\]
holds, where
\[
f_{\tilde V}=\langle\tilde \nabla^{\bot}\tilde V,\tilde\nabla^{\bot}\tilde V\rangle+\langle R(\tilde V),
\tilde V\rangle -\langle \tilde A^{\tilde V},\tilde A^{\tilde V}\rangle
\]
and $\tilde V$ is an arbitrary vector field of $\Gamma((T(L))^\bot)$ having compact
support on $L$ (see, for example, \cite{lw}).
\begin{theorem}
If $\calf$ is a minimal foliation of a manifold $M$ without boundary and the orthogonal distribution
$D^\bot$ is integrable, then any leaf $L$ of $\calf$ is stable.
\end{theorem}
{\it Proof.}
Let $\tilde V$ be an arbitrary vector field from $\Gamma((T(L))^\bot)$ having compact
support on $L$.
Since, for each point $q\in L$,  there exist a certain neighbourhood $\tilde U$ of $q$ in $L$ and
 $V\in\Gamma (D^\bot)$  such that $V|_{\tilde U}=\tilde V|_{\tilde U}$,
$W_{\tilde V}$ defined by
\begin{equation}
\label{e2}
W_{\tilde V}|_{\tilde U}=(\nabla_VV)^\top |_{\tilde U},
\end{equation}
is a well-defined  vector field of $\Gamma(TL)$.
Similarly, we can define\\ $\alpha_{\tilde V} \in L(\Gamma(TL),\Gamma((TL)^\bot))$ such that
\begin{equation}
\label{e3}
\alpha_{\tilde V}|_{\tilde U}=\alpha_V|_{\tilde U}.
\end{equation}
\par Now, let $p$ be  a fixed  point of $L$. Note that
\[
f_{\tilde V}(p)= f_{V,V}(p) ,
\]
where  $V$ as above. From Lemma \ref{l2}, we have
\[
f_{\tilde V}(p)= |\alpha_V|^2(p)-\divergence_L((\nabla_VV)^\top)(p).
\]
Using   (\ref{e2})  and (\ref{e3}),   we obtain
\[
f_{\tilde V}(p)= |\alpha_{\tilde V}|^2(p)-\divergence_L(W_{\tilde V})(p).
\]
Since  the point $p$ is  artbitrary, we get
\[
\int_L f_{\tilde V}=\int_L |\alpha_{\tilde V}|^2\geq 0.
\]
This ends the proof.\hfill$\square$
\par
Note that, the above theorem can be proved using the notion of  calibration \cite{ls}.
In our case, the volume form of leaves, which is a smooth $n$-form on $M$, gives a calibration
of $\calf$ (see \cite{hl}).
\section{Jacobi and Killing fields}
For the mapping $A: \Gamma(D^\bot )\rightarrow L(\Gamma(D),\Gamma(D))$ defined by
\[
A(V)=A^V,\quad V\in\Gamma(D^\bot),
\]
we can construct  $A^t$,  the transpose of $A$, i.e., if $B\in L(\Gamma(D),\Gamma(D))$ then
\[
\langle A^t(B),V\rangle(p)= \langle A^V,B\rangle(p), \quad  p\in M.
\]
We then set
\[
\hat A=A^t\circ A.
\]
Furthermore, if $V\in \Gamma(D^\bot)$, we construct a new cross-section $\nabla^{\bot^2} V$
in $D^\bot$ by setting
\begin{equation}
\label{e4}
\nabla^{\bot^2} V=\nabla_{e_i}^\bot\nabla^\bot_{e_i}V-\nabla^\bot_{\nabla_{e_i}^\top e_i} V,
\end{equation}
i.e., the trace of the connection of the mapping $\nabla^\bot V$.
\par Finally, we define $J:\Gamma(D^\bot)\rightarrow \Gamma(D^\bot)$ by
\[
J(V)=-\nabla^{\bot^2}V+R(V)-\hat A(V).
\]
\begin{definition}
{\rm We say that a normal section $V\in \Gamma(D^\bot)$  is a Jacobi field  of $\calf$ if $J(V)=0$ on $M$.}
\end{definition}
Similarly, we can  introduce  Jacobi fields of a leaf $L$ of the foliation $\calf$ (see \cite{s1}). Then
$V$ is a Jacobi  field of $\calf$ if, for any leaf $L$ of the foliation $\calf$,  $V|L$ is a  Jacobi field of $L$.
Moreover,  we will denote by $\alpha_V^t$  the transpose of $\alpha_V$, i.e.,
\[
\langle \alpha_V^t(W),X\rangle(p)= \langle W,\alpha_V (X)\rangle(p), \quad  p\in M
\]
for any  $X\in\Gamma(D)$, $W\in\Gamma(D^\bot)$.
\begin{lemma}
\label{l3}
Let $\calf$ be a minimal foliation of a manifold $M$ and assume that the orthogonal distribution is integrable. Then
we have the  formula
\[
\langle J(V),W\rangle=\langle \alpha_V,\alpha_W\rangle+\divergence_L(\alpha_V^t(W))
\]
for $V,W\in\Gamma(D^\bot)$.
\end{lemma}
{\it Proof.}
Let $p$ be a fixed point of $L$ and $(\{e_i\},\{ e_\alpha\})$ be a local adapted  frame field, such
that $\nabla_X^\top e_i(p)=0$ and $\nabla_X^\bot e_\alpha(p)=0$ for each vector  field $X$ on $M$.
Using Lemma  \ref{l2}  and the fact that  $\langle \hat A(V),W\rangle=\langle A^V,A^W\rangle$,
we obtain
\[
\langle J(V),W \rangle=-\langle \nabla^{\bot^2}V,W\rangle -\langle \nabla^\bot V,
\nabla^\bot W\rangle-\divergence_L((\nabla_VW)^\top)+\langle \alpha_V,\alpha_W\rangle.
\]
Then using  (\ref{e4}), we have the following equalities at the point $p$:
\begin{align*}
\langle J(V),W \rangle=&-\langle \nabla^\bot_{e_i}\nabla_{e_i}^\bot V,W\rangle -\langle \nabla^\bot V,
\nabla^\bot W\rangle-\divergence_L((\nabla_VW)^\top)+\langle \alpha_V,\alpha_W\rangle\\
=&-e_i\langle \nabla_{e_i}^\bot V,W\rangle-\divergence_L((\nabla_VW)^\top)+\langle \alpha_V,\alpha_W\rangle\\
=&-e_i\langle (\nabla_Ve_i)^\bot,W\rangle +e_i\langle \alpha_V(e_i), W\rangle\\
-&\divergence_L((\nabla_VW)^\top)+\langle \alpha_V,\alpha_W\rangle\\
=& e_i\langle (\nabla_V W)^\top,e_i\rangle +e_i\langle \alpha_V(e_i), W\rangle
-\divergence_L((\nabla_VW)^\top)+\langle \alpha_V,\alpha_W\rangle\\
=&e_i\langle \alpha_V(e_i), W\rangle+\langle \alpha_V,\alpha_W\rangle\\
=&e_i\langle \alpha^t_V(W),e_i\rangle+\langle \alpha_V,\alpha_W\rangle\\
=&\langle \nabla_{e_i}(\alpha^t_V(W)),e_i\rangle+\langle \alpha_V,\alpha_W\rangle\\
=&\divergence_L(\alpha_V^t(W)) +\langle \alpha_V,\alpha_W\rangle.
\end{align*}
Since the point $p$ is arbitrary, we complete the proof.\hfill$\square$
\begin{proposition}
\label{p1}
Let  $\calf$  be a minimal foliation of a manifold  $M$ with  the  integrable orthogonal distribution.
If a  vector field $X$ on $M$  is foliation preserving, i.e., maps leaves to leaves.
Then $X^\bot$ is a Jacobi field.
\end{proposition}
{\it Proof.} Since $X$ is foliation preserving, we have
\[
[X,F]\in\Gamma(D)\quad {\rm for}\quad F\in\Gamma(D).
\]
Consequently, $\alpha_{X^\bot}=0$ . From Lemma \ref{l3}, for $V=X^\bot$ and an arbitrary
$W\in \Gamma(D^\bot)$,  we have
\[
\langle J(V),W\rangle=0.
\]
Thus $X^\bot$ is a Jacobi field.\hfill$\square$
\par
Conversely, under an additional assumption, we have  the following proposition.
\begin{proposition}
\label{p2}
Let $\calf$  be a minimal foliation of a manifold  $M$ such that  all leaves are closed and the orthogonal distribution is
integrable. If $X\in\Gamma(TM)$ is a vector field such that $X^\bot$
 is a Jacobi field, then $X$ is foliation preserving.
\end{proposition}
{\it Proof.} It suffices to show  that $\alpha_V=0$  for $V=X^\bot$.
Since $V$ is a Jacobi field,  from Lemma \ref{l3},   for any leaf $L$ of $\calf$, we have
\[
0=\int_L\langle J(V),V\rangle =\int_L(|\alpha_V|^2+\divergence_L(\alpha_V^t(V)))=\int_L |\alpha_V|^2.
\]
Consequently, $\alpha_V=0$ on $M$.\hfill$\square$
\begin{corollary}
\label{c2}
Let $\calf$ be as in Proposition $\ref{p2}$. If $X$ is a Killing vector field on $M$, then $X$ is foliation preserving.
\end{corollary}
{\it Proof.} Since the normal component of  the Killing vector field is  a Jacobi field for each leaf $L$  (see \cite{s1}),  $X^\bot$ is a Jacobi field for each
$L$ and hence for $\calf$.\hfill$\square$
\par
Using the notion of  calibration, the above corollary was proved  by Oshikiri \cite{o3}.
\begin{remark}
Corollary \ref{c2} can not be extended to the case when the orthogonal
distribution is not integrable:
The Hopf fibration of the unit sphere $S^3 \rightarrow S^2$ gives a counter example.
\end{remark}
\begin{proposition}
\label{p3}
Let $\calf$ and $\calf^\bot$ be minimal orthogonal foliations on a closed manifold $M$.
If $X\in\Gamma(TM)$ is a vector field such that $X^\bot$
is a Jacobi field, then $X$ preserves $\calf$.
\end{proposition}
{\it Proof.} Denote by $H^\bot$ the mean curvature vector field of $\calf^\bot$, then, from
 Lemma \ref{l3} for
$V=W=X^\bot$,  we obtain
\begin{align*}
&0=\int_M\langle J(V),V\rangle=\int_M(|\alpha_V|^2+\divergence_L(\alpha_V^t(V)))\\
&=\int_M(|\alpha_V|^2+\divergence_M(\alpha_V^t(V))+\langle\alpha^t_V,H^\bot\rangle)\\
&=\int_M |\alpha_V|^2.
\end{align*}
Thus $\alpha_V=0$ and  $X$ preserves $\calf$.\hfill$\square$
\begin{corollary}
\label{c3}
Let $\calf,\calf^\bot$ be as in  Proposition $\ref{p3}$. If $X$ is a Killing vector field on $M$, then $X$ is foliation preserving.
\end{corollary}
\begin{proposition}
Let $\calf$ be a foliation  with all leaves compact of a   manifold  $M$ and $X\in\Gamma(TM)$
 a vector field preserving $\calf$.
If $X^{\bot}(p) =0$  and $p\in L\in \calf$, then $X^{\bot}=0$ on $L$.
\end{proposition}
{\it Proof.} Denote $V=X^\bot$, then $\alpha_V=0$ on $M$ and $V(p)=0$. Let $q$, $q \neq p$,
be  an arbitrary point of $L$.   Since $L$ is complete, there exists a geodesic $c:(-\epsilon,1+\epsilon)\rightarrow L$
 connecting $p$ and $q$ with $c(0)=p$ and $c(1)=q$ and a
covering  $\{U^I\}_{I=0}^{N}$ of $c$,  with an orthonormal frame $\{e_\alpha^I\}$.
Let $u^I=c^{-1}(U^I)$ for $I=0,\ldots,N$, and let $C^I\in \Gamma(D|_{U^I})$ be  vector fields such that $C^I(c(t))=\dot c(t)$ for $t\in u^I$.
From the assumption, for each $I$, we have
\[
\langle \alpha_V(C^I),e_\alpha^I\rangle (c(t))=0.
\]
Since $V=V_I^\alpha e^I_\alpha$ (summation over $\alpha$) on $U^I$, we have
\[
\frac{d}{dt}(V_I^\alpha \circ  c)(t)+(V_I^\beta\circ c)(t)(A^I)_{\beta}^\alpha(t)=0,
\]
for a matrix $A^I(t)$.
Thus $V_I^\alpha \circ  c$ is a solution of the set of  linear differential equations.
Without loss of generality, we may assume that  $p\in U^0$ and $q\in U^N$. Consequently,
 $V_0^\alpha \circ  c\equiv 0$. Inductively $V_{I}^\alpha \circ  c\equiv 0$ for any $I$,  and thus $V(q)=0$. \hfill$\square$
\begin{corollary}
Let $M$ be a  manifold and  $\calf$ a minimal foliation  having  all leaves
closed  and  the integrable  orthogonal distribution. If a Killing vector field $X$ on $M$ is tangent  to a leaf
$L$ at some point, then  $X$ is tangent  to $L$
everywhere on $L$.
\end{corollary}
\section*{Acknowledgment}
The author is grateful to  Pawe\l \,  Walczak and  Gen-ichi Oshikiri
for  disccusion   and  helpful comments.

Krzysztof Andrzejewski\\
Institute of Mathematics, Polish Academy of Sciences\\
ul. \'Sniadeckich 8, 00-956 Warszawa, Poland\\
{\it and}\\
Department of Theoretical Physics II,
University  of \L \' od\'z\\
ul. Pomorska 149/153, 90 - 236  \L \' od\'z, Poland.\\
{\it e-mail:} k-andrzejewski@uni.lodz.pl
\end{document}